\documentclass{amsart}

\usepackage{amsmath}
\usepackage{amsfonts}
\usepackage{amssymb}
\usepackage{graphicx}
\usepackage{hyperref}
\usepackage{mathrsfs}
\usepackage{pb-diagram}
\usepackage{epstopdf}
\usepackage{amsmath,amsfonts,amsthm,enumerate,amscd,latexsym}
\usepackage{curves}
\usepackage{bbm}
\usepackage[mathscr]{eucal}
\usepackage{caption}
\usepackage{subcaption}
\usepackage{tikz}
\newtheorem{thm}{Theorem}[section]

\newtheorem{prop}[thm]{Proposition}
\newtheorem{defn}[thm]{Definition}

\newtheorem{remark}[thm]{Remark}

\numberwithin{equation}{section}

\def\bZ{\mathbb{Z}}

\def\bR{\mathbb{R}}
\def\bC{\mathbb{C}}
\def\bP{\mathbb{P}}

\begin{document}

\title[HMS for local CY via SYZ]{Homological mirror symmetry\\ for local Calabi-Yau manifolds via SYZ}
\author[K. Chan]{Kwokwai Chan}
\address{Department of Mathematics\\ The Chinese University of Hong Kong\\ Shatin\\ Hong Kong}
\email{kwchan@math.cuhk.edu.hk}

\date{\today}

\begin{abstract}
This is a write-up of the author's talk in the conference {\em Algebraic Geometry in East Asia 2016} held at the University of Tokyo in January 2016. We give a survey on the series of papers \cite{Chan13, Chan-Ueda13, CPU16, CPU16a} where the author and his collaborators Daniel Pomerleano and Kazushi Ueda show how Strominger-Yau-Zaslow (SYZ) transforms can be applied to understand the geometry of Kontsevich's homological mirror symmetry (HMS) conjecture for certain local Calabi-Yau manifolds.
\end{abstract}

\maketitle

\tableofcontents

\section{Introduction}


The {\em homological mirror symmetry (HMS) conjecture}, proposed by Kontsevich in his 1994 ICM address \cite{Kontsevich-ICM94}, asserts that the derived Fukaya category of a Calabi-Yau manifold $X$ should be equivalent to the derived category of coherent sheaves over its mirror $\check{X}$ and vice versa. This gives a beautiful categorical formulation of mirror symmetry which is expected to imply the original enumerative predictions of mirror symmetry (see the recent work of Ganatra-Perutz-Sheridan \cite{GPS15}). The conjecture has been verified for elliptic curves \cite{Polishchuk-Zaslow98}, abelian varieties \cite{Fukaya02c, Kontsevich-Soibelman01, Abouzaid-Smith10}, quartic surfaces \cite{Seidel-K3}, quintic 3-folds \cite{Nohara-Ueda12} and higher-dimensional Calabi-Yau hypersurfaces in projective spaces \cite{Sheridan11b}.

There is also a formulation of the HMS conjecture for a Fano manifold $X$, again proposed by Kontsevich \cite{Kontsevich-ENS98}. In this case, the mirror is given by a Landau-Ginzburg (LG) model $(\check{X}, W)$ -- a pair consisting of a noncompact K\"ahler manifold $\check{X}$ and a holomorphic function $W:\check{X} \to \bC$ (called the superpotential). More generally, one can consider manifolds with an effective anti-canonical divisor \cite{Auroux07} or even general type manifolds \cite{Katzarkov07, KKOY09, GKR12}, for which mirrors are again given by LG models. In this setting, the HMS conjecture has been verified for
$\bP^2$ and $\bP^1\times\bP^1$ \cite{Seidel01a, Seidel01b}, toric del Pezzo surfaces \cite{Ueda06, Ueda-Yamazaki13}, weighted projective planes and Hirzebruch surfaces \cite{AKO08}, del Pezzo surfaces \cite{AKO06}, projective spaces \cite{Fang08} and more general projective toric varieties \cite{Abouzaid06, Abouzaid09, FLTZ11b, FLTZ12, FLTZ09},\footnote{A forthcoming work by Abouzaid-Fukaya-Oh-Ohta-Ono \cite{AFOOO} proves the other direction of the HMS conjecture for compact toric manifolds, building on previous works of Fukaya-Oh-Ohta-Ono \cite{FOOO-book, FOOO-toricI, FOOO-toricII, FOOO-toricIII}.} higher genus Riemann surfaces \cite{Seidel11, Efimov12} and Fano hypersurfaces in projective spaces \cite{Sheridan16}.\footnote{There are also many papers on the study of HMS for singularities which we do not attempt to list here.}

The proofs of these HMS statements, though often involve deep and ingenious arguments (e.g. the deformation technique pioneered by Seidel \cite{Seidel-ICM, Seidel-book}), are usually done by separate computations on the two sides and an explicit identification between the generators. In particular, there is often little discussion on explicit constructions of geometric functors implementing the HMS equivalences and the underlying geometry is hidden.

On the other hand, the proposal \cite{SYZ96} by Strominger, Yau and Zaslow in 1996 provides a geometric explanation for mirror symmetry -- the celebrated {\em SYZ conjecture} asserts that a pair of mirror manifolds should admit fiberwise dual Lagrangian torus fibrations (SYZ fibrations) over the same base and there should exist Fourier-Mukai--type transforms (SYZ transforms) responsible for the interchange of the symplectic and complex geometries over the mirror manifolds.

In the case when quantum corrections are absent, this has been exploited to understand the geometry of the HMS conjecture. Namely, for dual Lagrangian torus bundles over affine manifolds, a geometric functor realizing the HMS equivalences was constructed by Arinkin-Polishchuk \cite{Arinkin-Polishchuk01} and Leung-Yau-Zaslow \cite{LYZ00} using SYZ transforms (cf. the recent work by Tu \cite{Tu15} and Abouzaid \cite{Abouzaid14a, Abouzaid14b}). Their functor is a real version of the Fourier-Mukai transform on families of real tori and it takes a Lagrangian section to a holomorphic line bundle over the mirror. This was applied in the study of HMS for elliptic curves \cite{Arinkin-Polishchuk01}, abelian varieties \cite{Fukaya02c, Kontsevich-Soibelman01} and toric varieties \cite{Abouzaid06, Abouzaid09, Fang08, FLTZ11b, FLTZ12, FLTZ09, Tu15}.

To extend the constructions in \cite{Arinkin-Polishchuk01, LYZ00} to more general cases, a major difficulty lies in the presence of singularities in the SYZ fibrations. When singular fibers exist, the construction of the mirror manifold has to be modified by instanton corrections coming from holomorphic disks bounded by the regular Lagrangian fibers, as have been shown by Fukaya \cite{Fukaya05}, Kontsevich-Soibelman \cite{Kontsevich-Soibelman06}, Gross-Siebert \cite{Gross-Siebert-logI, Gross-Siebert-logII, Gross-Siebert-reconstruction, Gross-Siebert11b, Gross-Siebert12, Ruddat-Siebert14} and Auroux \cite{Auroux07, Auroux09}. Accordingly, the constructions of the desired functors have to be modified. In this article, we will demonstrate how this can be done when $X$ is the mirror of a toric Calabi-Yau manifold $\check{X}$, following \cite{Chan13, Chan-Ueda13, CPU16, CPU16a} (cf. also Gross-Matessi \cite{Gross-Matessi15}). Our slogan is:

\begin{center}
{\em ``The geometry of HMS is revealed by SYZ.''}
\end{center}

\section*{Acknowledgment}
The author would like to thank Daniel Pomerleano and Kazushi Ueda for very fruitful collaborations which lead to much of the work described here. He is also grateful to the organizers of the conference ``Algebraic Geometry in East Asia 2016'' held at the University of Tokyo in January 2016 for invitation and hospitality.

The research described in this article was substantially supported by a grant from the Research Grants Council of the Hong Kong Special Administrative Region, China (Project No. CUHK400213).

\section{Semi-flat SYZ transform}

To begin with, let us briefly recall the constructions by Arinkin-Polishchuk \cite{Arinkin-Polishchuk01} and Leung-Yau-Zaslow \cite{LYZ00}; the exposition here mainly follows the latter.

Let $B$ be an $n$-dimensional {\em tropical manifold}, meaning that the transition maps are elements of the group $\bR^n \rtimes \text{GL}_n(\bZ)$. Denote by $\Lambda\subset TB$ the bundle of lattices locally generated by the coordinate vector fields $\frac{\partial}{\partial x_1},\ldots,\frac{\partial}{\partial x_n}$, where $x_1,\ldots,x_n$ are (local) affine coordinates on $B$. Dually, let $\Lambda^\vee \subset T^*B$ be the family of lattices locally generated by the coordinate one-forms $dx_1,\ldots,dx_n$. Then we have a pair of manifolds
\begin{equation*}
X := T^*B/\Lambda^\vee,\quad \check{X} := TB/\Lambda,
\end{equation*}
together with fibrations
\begin{equation*}
\rho:X \to B, \quad \check{\rho}:\check{X} \to B,
\end{equation*}
which are fiberwise dual to each other.

Observe that $X$ is naturally a symplectic manifold: Let $\xi_1, \ldots, \xi_n$ be the fiber coordinates on $T^*B$, i.e. $(x_1, \ldots, x_n, \xi_1, \ldots, \xi_n)\in T^*B$ represents the cotangent vector $\sum_{j=1}^n \xi_j dx_j$ at the point $\mathbf{x} = (x_1, \ldots, x_n) \in B$. Then the canonical symplectic structure
\begin{equation*}
\omega := \sum_{j=1}^n dx_j\wedge d\xi_j
\end{equation*}
on $T^*B$ descends to a symplectic structure on $X$ (which, by abuse of notation, will still be denoted by $\omega$).

On the other hand, $\check{X}$ is naturally a complex manifold: Let $y_1, \ldots, y_n$ be the dual fiber coordinates on $TB$, i.e. $(x_1, \ldots, x_n, y_1, \ldots, y_n)\in TB$ represents the tangent vector $\sum_{j=1}^n y_j\frac{\partial}{\partial x_j}$ at the point $\mathbf{x} \in B$. Then the complex coordinates on $\check{X}$ are given by
\begin{equation*}
w_j := \exp 2\pi\left(x_j + \mathbf{i}y_j\right),\quad j = 1, \ldots, n.
\end{equation*}
Note that $\check{X}$ is equipped with the holomorphic volume form
\begin{equation*}
\check{\Omega} := d\log w_1 \wedge \cdots \wedge d\log w_n.
\end{equation*}

We can also equip $X$ and $\check{X}$ with compatible complex and symplectic structures respectively by choosing a {\em Hessian metric} (a metric given locally by the Hessian of a convex function) on $B$, hence giving both $X$ and $\check{X}$ a K\"ahler structure. If we further assume that the metric on $B$ satisfies the {\em real Monge-Amp\`ere equation}, we would obtain $T^n$-invariant Ricci-flat metrics on $X$ and $\check{X}$ making them {\em semi-flat Calabi-Yau manifolds}. Then the SYZ conjecture \cite{SYZ96} suggests that $X$ and $\check{X}$ form a mirror pair; see \cite{Leung05}, \cite{Chan-Leung10b} or \cite[Chapter 6]{D-branes-MS_book} for more discussion on this so-called {\em mirror symmetry without corrections}.

To construct the semi-flat SYZ transform \cite{Arinkin-Polishchuk01, LYZ00}, we first view the dual $T^*$ of a torus $T$ as the moduli space of flat $U(1)$-connections on the trivial line bundle $\underline{\bC}:=\bC\times T$ over $T$. If we are now given a section
$$L = \left\{ \left(\mathbf{x}, \mathbf{\xi}(\mathbf{x}) \right) \in X \mid \mathbf{x} \in B \right\}$$
of the fibration $\rho: X \to B$, then each point $\left(\mathbf{x}, \mathbf{\xi}(\mathbf{x})\right)\in L$ corresponds to a flat $U(1)$-connection $\nabla_{\mathbf{\xi}(\mathbf{x})}$ over the dual fiber $\check{\rho}^{-1}(\mathbf{x})\subset \check{X}$. The family of connections
$$\left\{ \nabla_{\mathbf{\xi}(\mathbf{x})} \mid \mathbf{x} \in B \right\}$$
patch together to give a $U(1)$-connection $\check{\nabla}$ over $\check{X}$. A straightforward calculation (which we leave to the reader) then shows the following

\begin{prop}
The $(0,2)$-part of the curvature two-form $F$ of $\check{\nabla}$ vanishes if and only if $L$ is Lagrangian with respect to $\omega$.
\end{prop}

Hence $\check{\nabla}$ determines a holomorphic line bundle $\mathcal{L}$ over $\check{X}$ precisely when $L$ is a Lagrangian section.

\begin{defn}\label{defn:semi-flat_SYZ_transform}
We define the {\em semi-flat SYZ transform} of the Lagrangian section $L$ in $X$ to be the holomorphic line bundle $\mathcal{L}$ over $\check{X}$, i.e.
\begin{equation*}
\mathcal{F}^\text{semi-flat}(L) := \mathcal{L}.
\end{equation*}
\end{defn}

\begin{remark}
More generally, one can consider Lagrangian submanifolds obtained from the {\em conormal bundle construction} and also equip them with flat $U(1)$-connections; see \cite[Sections 3.2 and 3.3]{LYZ00} or \cite[Section 2]{Chan13}.
\end{remark}


The semi-flat SYZ transform can be applied to understand the mirror symmetry and in particular the HMS conjecture for compact toric manifolds \cite{Abouzaid06, Abouzaid09, Chan-Leung10a, Chan-Leung10b, Fang08, FLTZ11b, FLTZ12, FLTZ09, Tu15}.

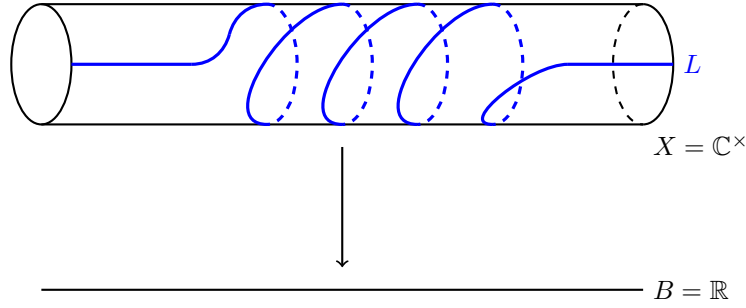
\begin{figure}
\begin{tikzpicture}
\draw [thick] (0,6) ellipse (0.4 and 0.8);
\draw [thick] (0,5.2) -- (8,5.2);
\draw [thick] (0,6.8) -- (8,6.8);
\draw [dashed, thick] (8,6.8) arc (90:270:0.4 and 0.8);
\draw [thick] (8,6.8) arc (90:-90:0.4 and 0.8) node[below right] {$X = \bC^\times$};
\draw [very thick, blue] (0.4,6) -- (2,6);
\draw [very thick, blue] (2,6) to [out=0, in=255] (2.5,6.4);
\draw [very thick, blue] (2.5,6.4) to [out=75, in=180] (3,6.8);
\draw [dashed, very thick, blue] (3,6.8) arc (90:-90:0.4 and 0.8);
\draw [very thick, blue] (3,5.2) to [out=180, in=180] (4,6.8);
\draw [dashed, very thick, blue] (4,6.8) arc (90:-90:0.4 and 0.8);
\draw [very thick, blue] (4,5.2) to [out=180, in=180] (5,6.8);
\draw [dashed, very thick, blue] (5,6.8) arc (90:-90:0.4 and 0.8);
\draw [very thick, blue] (5,5.2) to [out=180, in=180] (6,6.8);
\draw [dashed, very thick, blue] (6,6.8) arc (90:-90:0.4 and 0.8);
\draw [very thick, blue] (6,5.2) to [out=180, in=180] (7,6);
\draw [very thick, blue] (7,6) -- (8.4,6) node[right] {$L$};
\draw [thick, ->] (4,4.9) -- (4,3.3);
\draw [thick] (0,3) -- (8,3) node[right] {$B = \bR$};
\end{tikzpicture}
\caption{A Lagrangian section in $\bC^\times$ with prescribed asymptotic behavior.}\label{fig:Lag_Cstar}
\end{figure}

For example, let $B$ be the real line $\bR$ equipped with its natural $\bZ$-affine structure. Then $X = T*B/\Lambda^\vee$ is given by the punctured complex plane
$$T^*\bR/\bZ = \bC^\times = \bC \setminus \{0\}$$
with coordinates $x, \xi$ and symplectic structure $\omega = dx \wedge d\xi$, and the fibration $\rho: X \to B$ is nothing but the log function
$$\text{Log}: \bC^\times \to \bR, \quad z:= x + \mathbf{i}\xi \mapsto \log |z|.$$
In this case, a Lagrangian section $L$ (with certain asymptotic behavior controlled by the Laurent polynomial $W = z + \frac{q}{z}$) is a real curve like the one shown in Figure \ref{fig:Lag_Cstar}.

Its semi-flat SYZ transform produces a holomorphic line bundle
$$\mathcal{L} = \mathcal{F}^\text{semi-flat}(L)$$
over the toric variety $Y_\Sigma = \bP^1$; here the mirror
$$\check{X} = T\bR/\bZ \cong \bC^\times$$
should be regarded as the open dense orbit in $Y_\Sigma$. Furthermore, the degree of $\mathcal{L}$ is precisely given by the winding number $w(L)$ of $L$ (or its number of intersections with a radial axis, say, the negative real axis) so that
$$\mathcal{L} = \mathcal{O}_{\bP^1}(w(L)).$$

More generally, the mirror of a smooth projective toric variety $Y_\Sigma$ is given by an LG model $(X, W)$ where $X = (\bC^\times)^n$ and $W \in \bC[z_1^{\pm1}, \ldots, z_n^{\pm1}]$ is a Laurent polynomial whose Newton polytope coincides with the fan polytope (convex hull of the primitive generators of rays in $\Sigma$). In this case, (one direction of) the HMS conjecture asserts that the derived Fukaya-Seidel category $D^\pi Fuk(X, W)$ is equivalent to the derived category of coherent sheaves $D^b Coh(Y_\Sigma)$.

From the SYZ perspective, one wants to view $X = (\bC^\times)^n$ as $T^*B/\Lambda^\vee$ where $B = \bR^n$ is equipped with its natural $\bZ$-affine structure. The log map
$$\text{Log}: (\bC^\times)^n \to \bR^n, \quad (z_1, \ldots, z_n) \mapsto \left(\log|z_1|, \ldots, \log|z_n|\right)$$
provides a natural SYZ fibration (which is nothing but the projection $X = T^*B/\Lambda^\vee \to B$). The semi-flat SYZ dual
$$\check{X} = TB/\Lambda \cong (\bC^\times)^n$$
is then precisely the open dense torus orbit sitting inside $Y_\Sigma$. Now the semi-flat SYZ transform $\mathcal{F}^\text{semi-flat}$ again carries Lagrangian sections in $X$, with certain asymptotic behavior specified by the superpotential $W$, to holomorphic line bundles over $\check{X}$ which can be extended over $Y_\Sigma$ \cite{Chan09}.

\begin{remark}
For the applications of SYZ transforms to the other direction of HMS for compact toric manifolds, we refer the readers to \cite{Chan-Leung12, CHL12, CHL14}.
\end{remark}

\section{SYZ transforms for local Calabi-Yau manifolds}

To see how SYZ transforms can be constructed when quantum corrections are present, we will consider mirror symmetry for certain local Calabi-Yau manifolds.\footnote{In this article, a variety $Y$ is called {\em Calabi-Yau} if it is Gorenstein and its canonical line bundle $K_Y$ is trivial.}

\subsection{Local mirror symmetry by SYZ}

Let
$$f = f(z_1, \ldots, z_{n-1}) \in \bC[z_1^{\pm1}, \ldots, z_{n-1}^{\pm1}]$$
be a Laurent polynomial in $n-1$ variables with complex coefficients. Then the hypersurface
\begin{equation*}
X = \left\{ (x, y, z_1, \ldots, z_{n-1}) \in \bC^2 \times (\bC^\times)^{n-1} \mid xy = f(z_1, \ldots, z_{n-1}) \right\}
\end{equation*}
in $\bC^2 \times (\bC^\times)^{n-1}$
is a noncompact Calabi-Yau variety since it admits the following holomorphic volume form:
$$\Omega := \text{Res}\left[ \frac{dx \wedge dy \wedge d\log z_1 \wedge \cdots \wedge d\log z_{n-1}}{xy - f(z_1,\ldots,z_{n-1})} \right].$$

It has long been known that the mirror of $X$ should be given by a {\em toric Calabi-Yau manifold}.
Let $N \cong \bZ^n$ be a rank $n$ lattice. Let $Y = Y_\Sigma$ be a toric variety defined by a fan $\Sigma$ in $N_\bR = N \otimes_\bZ \bR$.
The primitive generators of rays in $\Sigma$ will be denoted by $\nu_1, \ldots,\nu_m \in N$ which correspond to the toric prime divisors $D_1, \ldots, D_m \subset Y$ respectively.
The toric variety $Y$ is Calabi-Yau if and only if there exists a lattice point $\eta \in M = \text{Hom}(N, \bZ)$ such that $\langle \eta, \nu_i\rangle = 1$ for $i = 1, \ldots, m$ \cite{CLS_toric_book}. This in turn is equivalent to the existence of $\eta \in M$ such that the corresponding character $\chi^\eta \in \text{Hom}(M\otimes_\bZ \bC^\times, \bC^\times)$ defines a holomorphic function on $Y$ which has simple zeros along each toric prime divisor $D_i$ and is non-vanishing elsewhere. Note that $Y$ is necessarily noncompact in this case.

Let $Y = Y_\Sigma$ be a toric Calabi-Yau variety of complex dimension $n$.
We further assume that $Y$ is smooth and the fan $\Sigma$ has convex support. These conditions are satisfied if and only if $Y$ is a crepant resolution of an affine toric variety (defined by the cone $|\Sigma|$) with Gorenstein canonical singularities, which in turn is equivalent to saying that $Y$ is {\em semi-projective} \cite[p.332]{CLS_toric_book}.

Important examples of toric Calabi-Yau manifolds are given by total spaces of canonical line bundles $K_Z$ over compact toric manifolds $Z$. (In this case, semi-projectivity is equivalent to requiring that the base manifold $Z$ is semi-Fano, i.e. its canonical line bundle $K_Z$ is nef.) For example, the total space of $K_{\bP^2}=\mathcal{O}_{\bP^2}(-3)$ is a toric Calabi-Yau 3-fold whose fan $\Sigma$ has rays generated by the vectors
$$\nu_0 = (0, 0, 1),\quad \nu_1 = (1, 0, 1),\quad \nu_2 = (0, 1, 1),\quad \nu_3 = (-1, -1, 1) \in \bZ^3.$$

Mirror symmetry in this setting is called {\em local mirror symmetry} because it originated from applying mirror symmetry techniques to Fano surfaces (e.g $\bP^2$) contained inside compact Calabi-Yau manifolds, and could be derived using physical arguments from mirror symmetry for compact Calabi-Yau hypersurfaces in toric varieties by taking certain limits in the complexified K\"ahler and complex moduli spaces \cite{KKV97}. There is a large body of work, by both mathematicians and physicists, on this mirror symmetry; see e.g. \cite[Section 4]{Chan14} and the references therein.

In order to construct the mirror of $X$ using SYZ, one should find a Lagrangian torus fibration $\rho: X \to B$ with a Lagrangian section, and proceed in the following steps:\footnote{Such a procedure was pioneered by Auroux in \cite{Auroux07, Auroux09}, and later generalized to toric Calabi-Yau manifolds in \cite{CLL12} and certain blowups of toric varieties in \cite{AAK12}; see also \cite{Gross-Siebert_ICM}.}
\begin{itemize}
\item[Step 1]
Over the smooth locus $B_0 \subset B$, the pre-image
$$X_0 := \rho^{-1}(B_0)$$
can be identified with the quotient $T^*B_0/\Lambda^\vee$ by Duistermaat's action-angle coordinates \cite{Duistermaat80}.

\item[Step 2]
Then we have the {\em semi-flat} mirror
$$\check{X}_0 := TB_0/\Lambda,$$
which is naturally a complex manifold because the fibration $\rho$ induces a tropical affine structure on $B_0$.
However, this is not quite the correct mirror because we want to compactify $\check{X}_0$ but the complex structure on $\check{X}_0$ {\em cannot} be extended further to {\em any} (even partial) compactification due to nontrivial monodromy of the tropical affine structure around the discriminant locus $\Gamma = B \setminus B_0$.

\item[Step 3]
To obtain the correct and (partially) compactified mirror $\check{X}$, we need to modify the complex structure on $\check{X}_0$ by {\em instanton corrections} coming from holomorphic disks in $X$ bounded by the Lagrangian torus fibers of $\rho$.
\end{itemize}

For the examples we consider here, such a construction was carried out by Abouzaid-Auroux-Katzarkov \cite{AAK12}. The resulting {\em SYZ mirror} of $X$ is given by
$$\check{X} = Y_\Sigma \setminus H,$$
where $Y_\Sigma$ is a toric Calabi-Yau $n$-fold whose fan consists of cones over a regular subdivision of the Newton polytope of $f$, and the hypersurface $H$ is a smoothing of the union of toric prime divisors in $Y_\Sigma$.

For example, the mirror of
\begin{equation*}
X = \left\{ (x, y, z_1, z_2) \in \bC^2 \times (\bC^\times)^2 \mid xy = 1 + z_1 + z_2 \right\}
\end{equation*}
is given by
$$\bC^3 \setminus \{w_1w_2w_3 = 1\},$$
where $w_1, w_2, w_3$ are standard coordinates on $(\bC)^3$), while the mirror of
\begin{equation*}
X = \left\{ (x, y, z_1, z_2) \in \bC^2 \times (\bC^\times)^2 \mid xy = t + z_1 + z_2 + \frac{1}{z_1z_2} \right\}
\end{equation*}
is given by
$$K_{\bP^2} \setminus H,$$
where $H$ is a smoothing of the union of the four toric prime divisors in $K_{\bP^2}$.
We will briefly review the construction in \cite{AAK12} in some explicit examples below.

\subsection{The 2-dimensional case}

Let us consider
\begin{equation*}
X = \left\{ (x, y, z) \in \bC^2 \times \bC^\times \mid xy = f(z) \right\},
\end{equation*}
where $f(z) \in \bC[z, z^{-1}]$ is a Laurent polynomial in one variable. We equip $X$ with the symplectic structure given by restriction of the standard one
\begin{equation*}
\omega = -\frac{\mathbf{i}}{2}\left(dx\wedge d\bar{x} + dy\wedge d\bar{y} + \frac{dz\wedge d\bar{z}}{|z|^2}\right)
\end{equation*}
on $\bC^2 \times \bC^\times$ to $X$.

Then a Lagrangian torus fibration on $X$ is given explicitly by the map
\begin{align*}
\rho: X \to B := \bR^2,\quad (x, y, z) \mapsto \left( \log|z|, \mu(x, y, z) \right),
\end{align*}
where
\begin{equation*}
\mu(x, y, z)=\frac{1}{2}\left( |x|^2 - |y|^2 \right): X \to \bR
\end{equation*}
is the moment map associated to the Hamiltonian $S^1$-action:
\begin{equation*}
e^{\mathbf{i} \theta} \cdot (x, y, z) = (e^{\mathbf{i}\theta}x, e^{-\mathbf{i}\theta}y, z).
\end{equation*}

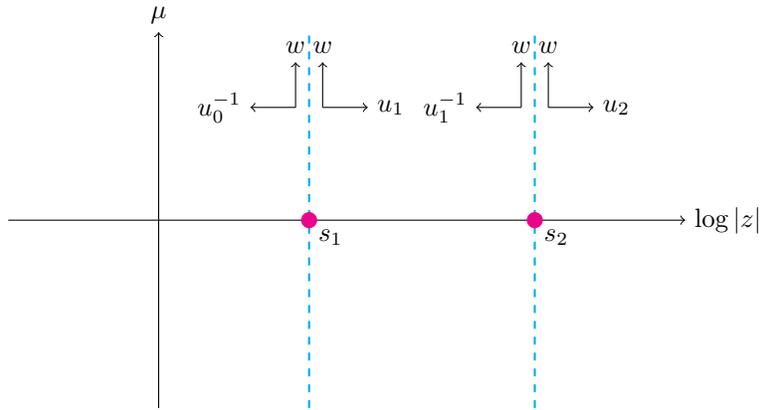
\begin{figure}
\begin{tikzpicture}
\draw[->] (1,3.5) -- (10,3.5) node[right] {$\log |z|$};
\draw[->] (3,1) -- (3,6) node[above] {$\mu$};
\draw [dashed, thick, cyan] (5,1) -- (5,6);
\draw [dashed, thick, cyan] (8,1) -- (8,6);
\draw[fill, magenta] (5,3.5) circle [radius=0.1] node[below right, black] {$s_1$};
\draw[fill, magenta] (8,3.5) circle [radius=0.1] node[below right, black] {$s_2$};
\draw[->] (4.82,5) -- (4.22,5) node[left] {$u_0^{-1}$};
\draw[->] (4.82,5) -- (4.82,5.6) node[above] {$w$};
\draw[->] (5.18,5) -- (5.78,5) node[right] {$u_1$};
\draw[->] (5.18,5) -- (5.18,5.6) node[above] {$w$};
\draw[->] (7.82,5) -- (7.22,5) node[left] {$u_1^{-1}$};
\draw[->] (7.82,5) -- (7.82,5.6) node[above] {$w$};
\draw[->] (8.18,5) -- (8.78,5) node[right] {$u_2$};
\draw[->] (8.18,5) -- (8.18,5.6) node[above] {$w$};
\end{tikzpicture}
\caption{The base of SYZ fibration for $X$ and its chamber structure}\label{fig:2d_SYZbase}
\end{figure}

For a more concrete example, we may consider
$$f(z) = (z - a)(z - b),$$
where $1 < a < b$ are real numbers. Then the base of the SYZ fibration is displayed in Figure \ref{fig:2d_SYZbase}.
Let $s_1 := \log a$ and $s_2 := \log b$. Then we have precisely two singular fibers (which are pinched 2-tori) over the points
\begin{equation*}
\Gamma = \{(s_1, 0), (s_2, 0)\} \subset B.
\end{equation*}
Also, the locus of Lagrangian torus fibers which bound nontrivial holomorphic disks is given by the union of two vertical lines (the cyan dashed lines in Figure \ref{fig:2d_SYZbase}):
\begin{equation*}
\{s_1\} \times \bR,\quad \{s_2\} \times \bR \subset B.
\end{equation*}
We call each of these vertical lines a {\em wall} in $B$. (For a more general $f(z)$ which has $k$ zeros (with distinct absolute values), there will be $k$ singular fibers and $k$ walls which are all vertical lines; all constructions which follow work in exactly the same way as in this particular example.)

\begin{figure}
\begin{tikzpicture}
\path [fill=lightgray] (-4,1) rectangle (-2,6);
\draw [dashed, thick] (-4,3.5) -- (-2,3.5);
\draw [dashed, thick, cyan] (-2,1) -- (-2,6);
\draw [fill, white] (-2,3.5) circle [radius=0.1];
\draw [thick, magenta] (-2,3.5) circle [radius=0.1];
\draw (-3,3.5) node[below] {$U_0$};
\path [fill=lightgray] (-1,1) rectangle (2,6);
\draw [dashed, thick] (-1,3.5) -- (2,3.5);
\draw [dashed, thick, cyan] (-1,1) -- (-1,6);
\draw [dashed, thick, cyan] (2,1) -- (2,6);
\draw [fill, white] (2,3.5) circle [radius=0.1];
\draw [thick, magenta] (2,3.5) circle [radius=0.1];
\draw [fill, white] (-1,3.5) circle [radius=0.1];
\draw [thick, magenta] (-1,3.5) circle [radius=0.1];
\draw (0.5,3.5) node[below] {$U_1$};
\path [fill=lightgray] (3,1) rectangle (5,6);
\draw [dashed, thick] (3,3.5) -- (5,3.5);
\draw [dashed, thick, cyan] (3,1) -- (3,6);
\draw [fill, white] (3,3.5) circle [radius=0.1];
\draw [thick, magenta] (3,3.5) circle [radius=0.1];
\draw (4,3.5) node[below] {$U_2$};
\end{tikzpicture}
\caption{Connected components of the complement of walls in $B$}\label{fig:2d_complement_wall}
\end{figure}
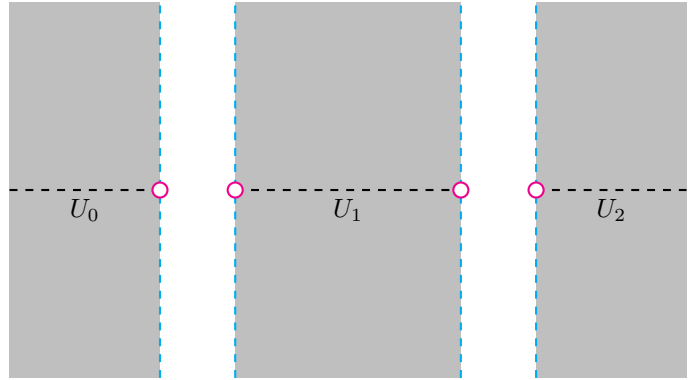

To see how the SYZ mirror is constructed, we consider the connected components (chambers) of the complement of the walls in $B$ (see Figure \ref{fig:2d_complement_wall}):
\begin{equation*}
U_0 := (-\infty, s_1) \times \bR,\quad U_1 := (s_1, s_2) \times \bR,\quad U_2 := (s_2, \infty) \times \bR.
\end{equation*}
Let $u_i$ ($i = 1, 2, 3, 4$) and $w$ be exponentiation of complexification of the affine coordinates as shown in Figure \ref{fig:2d_SYZbase}. Note that $w$ is a global coordinate since it comes from the moment map of the global Hamiltonian $S^1$-action on $X$.

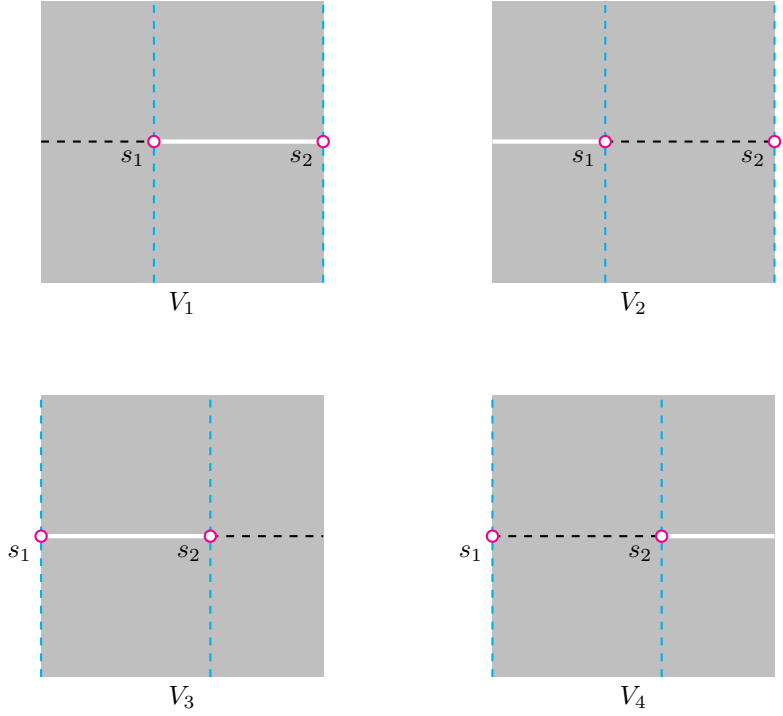
\begin{figure}
\begin{tikzpicture}
\path [fill=lightgray] (-3,0.75) rectangle (0.75,4.5);
\draw [ultra thick, white] (-1.5,2.625) -- (0.75,2.625);
\draw [dashed, thick] (-3,2.625) -- (-1.5,2.625);
\draw [dashed, thick, cyan] (-1.5,0.75) -- (-1.5,4.5);
\draw [dashed, thick, cyan] (0.75,0.75) -- (0.75,4.5);
\draw [fill, white] (-1.5,2.625) circle [radius=0.075] node[below left, black] {$s_1$};
\draw [thick, magenta] (-1.5,2.625) circle [radius=0.075];
\draw [fill, white] (0.75,2.625) circle [radius=0.075] node[below left, black] {$s_2$};
\draw [thick, magenta] (0.75,2.625) circle [radius=0.075];
\draw (-1.125,0.75) node[below] {$V_1$};
\path [fill=lightgray] (3,0.75) rectangle (6.75,4.5);
\draw [ultra thick, white] (3,2.625) -- (4.5,2.625);
\draw [dashed, thick] (4.5,2.625) -- (6.75,2.625);
\draw [dashed, thick, cyan] (4.5,0.75) -- (4.5,4.5);
\draw [dashed, thick, cyan] (6.75,0.75) -- (6.75,4.5);
\draw [fill, white] (4.5,2.625) circle [radius=0.075] node[below left, black] {$s_1$};
\draw [thick, magenta] (4.5,2.625) circle [radius=0.075];
\draw [fill, white] (6.75,2.625) circle [radius=0.075] node[below left, black] {$s_2$};
\draw [thick, magenta] (6.75,2.625) circle [radius=0.075];
\draw (4.875,0.75) node[below] {$V_2$};
\path [fill=lightgray] (-3,-4.5) rectangle (0.75,-0.75);
\draw [ultra thick, white] (-3,-2.625) -- (-0.75,-2.625);
\draw [dashed, thick] (-0.75,-2.625) -- (0.75,-2.625);
\draw [dashed, thick, cyan] (-3,-4.5) -- (-3,-0.75);
\draw [dashed, thick, cyan] (-0.75,-4.5) -- (-0.75,-0.75);
\draw [fill, white] (-3,-2.625) circle [radius=0.075] node[below left, black] {$s_1$};
\draw [thick, magenta] (-3,-2.625) circle [radius=0.075];
\draw [fill, white] (-0.75,-2.625) circle [radius=0.075] node[below left, black] {$s_2$};
\draw [thick, magenta] (-0.75,-2.625) circle [radius=0.075];
\draw (-1.125,-4.5) node[below] {$V_3$};
\path [fill=lightgray] (3,-4.5) rectangle (6.75,-0.75);
\draw [ultra thick, white] (5.25,-2.625) -- (6.75,-2.625);
\draw [dashed, thick] (3,-2.625) -- (5.25,-2.625);
\draw [dashed, thick, cyan] (3,-4.5) -- (3,-0.75);
\draw [dashed, thick, cyan] (5.25,-4.5) -- (5.25,-0.75);
\draw [fill, white] (3,-2.625) circle [radius=0.075] node[below left, black] {$s_1$};
\draw [thick, magenta] (3,-2.625) circle [radius=0.075];
\draw [fill, white] (5.25,-2.625) circle [radius=0.075] node[below left, black] {$s_2$};
\draw [thick, magenta] (5.25,-2.625) circle [radius=0.075];
\draw (4.875,-4.5) node[below] {$V_4$};
\end{tikzpicture}
\caption{The open sets $V_1$, $V_2$, $V_3$ and $V_4$ in $B$.}\label{fig:2d_open_sets_SYZbase}
\end{figure}

\begin{remark}
The reader may complain that the open sets $U_0$, $U_1$ and $U_2$ do not form an open cover of the smooth locus $B_0 = B \setminus \Gamma$. We are indeed cheating a little bit here because we want to simplify the exposition. To be more precise, one should consider, for example, the following open subsets in $B$ (see Figure \ref{fig:2d_open_sets_SYZbase}):
\begin{align*}
V_1 & := \left( (-\infty, s_2) \times \bR \right) \setminus \left( (s_1, s_2) \times \{0\} \right),\\
V_2 & := \left( (-\infty, s_2) \times \bR \right) \setminus \left( (-\infty, s_1) \times \{0\} \right),\\
V_3 & := \left( (s_1, \infty) \times \bR \right) \setminus \left( (s_1, s_2) \times \{0\} \right),\\
V_4 & := \left( (s_1, \infty) \times \bR \right) \setminus \left( (s_2, \infty) \times \{0\} \right).
\end{align*}
The following constructions work in exactly the same way using this set of open charts.
\end{remark}

Using the fact that the monodromy of the tropical affine structure going counter-clockwise around both $(s_1, 0)$ and $(s_2, 0)$ is given by the matrix
$$\left(
\begin{array}{ll}
1 & 1\\
0 & 1
\end{array}
\right),$$
we see that the semi-flat SYZ mirror $\check{X}_0$ is given by
$$\left(TU_0/\Lambda\right) \cup \left(TU_1/\Lambda\right) \cup \left(TU_2/\Lambda\right)$$
via the na\"ive gluing:
\begin{align*}
u_i & = u_{i+1} \text{ on $U_i \cap U_{i+1} \cap B_+$,}\\
u_i & = u_{i+1} w \text{ on $U_i \cap U_{i+1} \cap B_-$}
\end{align*}
for $i = 0, 1$, where $B_+$ (resp. $B_-$) is the upper half plane $\bR \times \bR_{>0}$ (resp. lower half plane $\bR \times \bR_{<0}$).

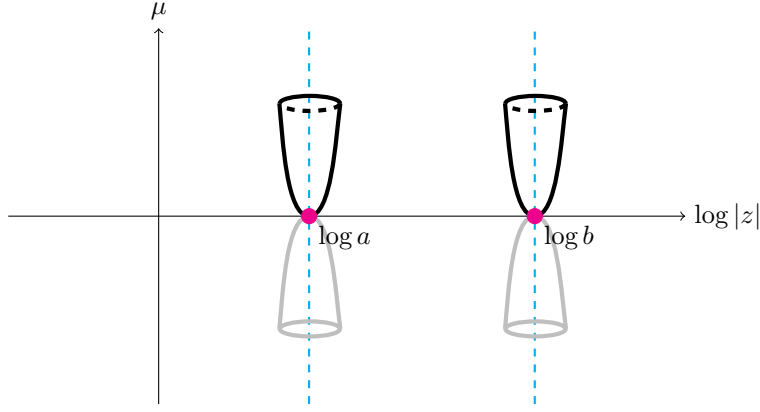
\begin{figure}
\begin{tikzpicture}
\draw [->] (1,3.5) -- (10,3.5) node[right] {$\log |z|$};
\draw [->] (3,1) -- (3,6) node[above] {$\mu$};
\draw [dashed, thick, cyan] (5,1) -- (5,6);
\draw [dashed, ultra thick] (5.41,5) arc (0:-180:0.4 and 0.1);
\draw [ultra thick] (5.41,5) arc (0:180:0.4 and 0.1);
\draw [ultra thick] (4.605,5) .. controls (4.65,4.5) and (4.7,3.5) .. (5,3.5);
\draw [ultra thick] (5.415,5) .. controls (5.35,4.5) and (5.3,3.5) .. (5,3.5);
\draw [ultra thick, lightgray] (5.41,2) arc (0:-180:0.4 and 0.1);
\draw [ultra thick, lightgray] (5.41,2) arc (0:180:0.4 and 0.1);
\draw [ultra thick, lightgray] (4.605,2) .. controls (4.65,2.5) and (4.7,3.5) .. (5,3.5);
\draw [ultra thick, lightgray] (5.415,2) .. controls (5.35,2.5) and (5.3,3.5) .. (5,3.5);
\draw [dashed, thick, cyan] (8,1) -- (8,6);
\draw [dashed, ultra thick] (8.41,5) arc (0:-180:0.4 and 0.1);
\draw [ultra thick] (8.41,5) arc (0:180:0.4 and 0.1);
\draw [ultra thick] (7.605,5) .. controls (7.65,4.5) and (7.7,3.5) .. (8,3.5);
\draw [ultra thick] (8.415,5) .. controls (8.35,4.5) and (8.3,3.5) .. (8,3.5);
\draw [ultra thick, lightgray] (8.41,2) arc (0:-180:0.4 and 0.1);
\draw [ultra thick, lightgray] (8.41,2) arc (0:180:0.4 and 0.1);
\draw [ultra thick, lightgray] (7.605,2) .. controls (7.65,2.5) and (7.7,3.5) .. (8,3.5);
\draw [ultra thick, lightgray] (8.415,2) .. controls (8.35,2.5) and (8.3,3.5) .. (8,3.5);
\draw [fill, magenta] (5,3.5) circle [radius=0.1] node[below right, black] {$\log a$};
\draw [fill, magenta] (8,3.5) circle [radius=0.1] node[below right, black] {$\log b$};
\end{tikzpicture}
\caption{Holomorphic disks bounded by fibers over the walls}\label{fig:2d_holom_disks}
\end{figure}

As we explained before, this is not quite the correct mirror because nontrivial monodromy around the discriminant locus $\Gamma$ makes it impossible to extend the complex structure of $\check{X}_0$ any further. The corrections we need come precisely from the nontrivial holomorphic disks bounded by fibers over the walls \cite{Auroux07, Auroux09, CLL12, AAK12} as shown in Figure \ref{fig:2d_holom_disks}. This leads to the corrected gluing given by:
\begin{equation}\label{eqn:corrected_gluing}
\begin{split}
u_i & = u_{i+1} (1 + w)  \text{ on $U_i \cap U_{i+1} \cap B_+$,}\\
u_i & = u_{i+1} w (1 + w^{-1})  \text{ on $U_i \cap U_{i+1} \cap B_-$},
\end{split}
\end{equation}
which can simply be written as
\begin{equation*}
u_i = u_{i+1} (1 + w)\text{ on $U_i \cap U_{i+1}$}
\end{equation*}
for $i = 0, 1$.

\begin{figure}
\begin{tikzpicture}
\draw [thick] (-2,4.5) -- (0,2);
\draw [thick] (0,2) -- (4,2);
\draw [thick] (4,2) -- (6,4.5);
\draw [dashed, thick, cyan] (0,2) -- (0,4.5);
\draw [dashed, thick, cyan] (4,2) -- (4,4.5);
\draw[->] (-0.2,2.5) -- (-0.7,3.125) node[above] {$u_0^{-1}$};
\draw[->] (-0.2,2.5) -- (-0.2,3.125) node[above] {$w$};
\draw[->] (0.2,2.5) -- (0.825,2.5) node[right] {$u_1$};
\draw[->] (0.2,2.5) -- (0.2,3.125) node[above] {$w$};
\draw[->] (3.8,2.5) -- (3.175,2.5) node[left] {$u_1^{-1}$};
\draw[->] (3.8,2.5) -- (3.8,3.125) node[above] {$w$};
\draw[->] (4.2,2.5) -- (4.7,3.125) node[above] {$u_2$};
\draw[->] (4.2,2.5) -- (4.2,3.125) node[above] {$w$};
\end{tikzpicture}
\caption{The corrected gluing in the SYZ mirror of $X$}\label{fig:2d_twisted_gluing}
\end{figure}
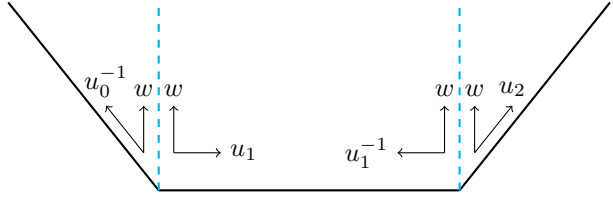

In fact this is exactly the gluing of complex charts in $T^*\bP^1$ (see Figure \ref{fig:2d_twisted_gluing}), which is the toric resolution of the $A_1$-singularity $\bC^2/\bZ_2$. Hence we conclude that the SYZ mirror of
\begin{equation*}
X = \left\{ (x, y, w) \in \bC^2 \times \bC^\times \mid xy = (w - a)(w - b) \right\}
\end{equation*}
is given by
\begin{equation*}
\check{X} = T^*\bP^1 \setminus H,
\end{equation*}
where $H$ is a smoothing of the union of the three toric divisors in $T^*\bP^1$.

\begin{figure}
\begin{tikzpicture}
\draw[->] (-0.75,2.625) -- (7.5,2.625);
\draw[->] (3,-0.375) -- (3,6);
\draw [very thick] (4.5,2.475) -- (4.5,2.775);
\draw (4.5,2.625) node[below right] {$a$};
\draw [very thick] (7.125,2.475) -- (7.125,2.775);
\draw (7.125,2.625) node[below right] {$b$};
\draw [very thick] (3,2.625) circle [radius=0.075];
\draw [ultra thick, blue] (3,2.625) to [out=0, in=180] (4.5,3.375);
\draw [ultra thick, blue] (4.5,3.375) to [out=0, in=90] (5.25,2.625);
\draw [ultra thick, blue] (5.25,2.625) to [out=270, in=0] (3.75,1.3125);
\draw [ultra thick, blue] (3.75,1.3125) to [out=180, in=270] (0.75,3);
\draw [ultra thick, blue] (0.75,3) to [out=90, in=180] (3.75,5.25);
\draw [ultra thick, blue] (3.75,5.25) to [out=0, in=90] (6.75,2.625);
\draw [ultra thick, blue] (6.75,2.625) to [out=270, in=0] (3.75,0);
\draw [ultra thick, blue] (3.75,0) to [out=180, in=315] (0,1.125);
\draw [ultra thick, blue] (0,1.125) to [out=135, in=0] (-1.875,2.625);
\end{tikzpicture}
\caption{An admissible path in $\bC^\times$}\label{fig:2d_admissible_path}
\end{figure}
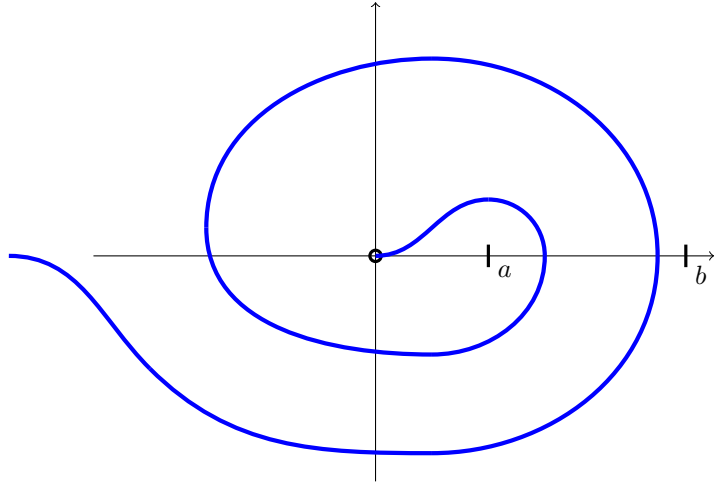

Now we turn to the construction of SYZ transforms in this example. First of all, Lagrangian sections of the SYZ fibration $\rho: X \to B$ can be constructed by considering (admissible) paths $\gamma$ in the base of the conic fibration
\begin{equation*}
\pi: X \to \bC^\times,\quad (u, v, z) \mapsto z.
\end{equation*}
An example is shown in Figure \ref{fig:2d_admissible_path}. Taking the subset in $X$ swept by the parallel transports of the real locus of a fixed conic fiber, we get a Lagrangian section $L = L_\gamma$ of $\rho: X \to B$ which is topologically a copy of $\bR^2$.

To construct the SYZ transform of $L$, we first note that the semi-flat SYZ transform of $L$ over each chamber $U_i$ ($i = 0, 1, 2$) gives a $U(1)$-connection $\check{\nabla}_i$ over the open chart $U_i^\bC := TU_i/\Lambda \cong (\bC^\times)^2$ in $\check{X}$, but they do not quite agree over the walls due to correction of the gluing \eqref{eqn:corrected_gluing}. More precisely, we have
$$\check{\nabla}_i - \check{\nabla}_{i-1} = - \mathbf{i} \xi_1(s_i) d\text{arg}(1 + w)$$
over the wall $\{s_i\} \times \bR$ for $i = 1, 2$.

The key observation is that, although the connections do not match over the walls, they only {\em differ by a gauge}! So after modification by adding a gauge, the connections $\check{\nabla}_0$ and $\check{\nabla}_1$ (resp. $\check{\nabla}_1$ and $\check{\nabla}_2$) can patch together to give us a $U(1)$-connection $\tilde{\nabla}_1$ (resp. $\tilde{\nabla}_2$) over the open chart $U_{01}^\bC := V_1^\bC \cup V_2^\bC$ (resp. $U_{12}^\bC := V_3^\bC \cup V_4^\bC$)).

This produces a global $U(1)$-connection over the mirror $\check{X}$, whose curvature has vanishing $(0,2)$-part again because $L$ is Lagrangian. Hence we can define the SYZ transform of $L$ as a holomorphic line bundle $\mathcal{L}$ over $\check{X}$. Moreover, in the overlap $U_{01}^\bC \cap U_{12}^\bC$, we have
$$\tilde{\nabla}_2 - \tilde{\nabla}_1 = \mathbf{i} (\xi(s_2) - \xi(s_1)) d\text{arg}(1 + w),$$
implying that the transition function for $\mathcal{L}$ from $U_{01}$ to $U_{02}$ is given by
$$(1 + w)^{-(\xi(s_2) - \xi(s_1))}.$$
Using this, one can show that the degree of $\mathcal{L}$ over the exceptional $\bP^1 \subset \check{X}$ (which determines the isomorphism class of the line bundle) is precisely given by the intersection number $\gamma \cdot [a, b]$.

Similar techniques work for all mirrors of $A_n$-resolutions (or toric Calabi-Yau 2-folds) \cite{Chan-Ueda13} and also for the deformed conifold and mirrors of other small toric Calabi-Yau 3-folds \cite{CPU16}. Moreover, this SYZ transform, which is defined only on the object level, in fact induces HMS equivalences as predicted by Kontsevich; the readers may consult \cite{Chan-Ueda13, CPU16} for the precise statements and proofs.

\subsection{The 3-dimensional case}

In the general 3-dimensional case, construction of Lagrangian torus fibrations becomes more subtle. Given
\begin{equation*}
X = \left\{ (x, y, z_1, z_2) \in \bC^2 \times (\bC^\times)^2 \mid xy = f(z_1, z_2) \right\},
\end{equation*}
where $f(z_1,z_2) \in \bC[z_1^{\pm1}, z_2^{\pm1}]$ is a Laurent polynomial in two variables,
we first regard $X$ as an open dense subset in the blowup
$$p: \text{Bl}_{0 \times Z}\left(\bC \times (\bC^\times)^2\right) \to \bC \times (\bC^\times)^2$$
of $\bC \times (\bC^\times)^2$ along the codimension two locus $0 \times Z$; here, $Z \subset (\bC^\times)^2$ is the hypersurface defined by $f(z_1,z_2)$, i.e.
$$Z = \{ (z_1, z_2) \in (\bC^\times)^2 \mid f(z_1, z_2) = 0\}.$$
Then we equip $X$ with the blowup symplectic structure $\omega_\epsilon$ where $\epsilon > 0$, the symplectic area of the exceptional $\bP^1$'s (i.e. fibers of $p$ over $0 \times Z$), is sufficiently small.\footnote{We do not use the restriction of the standard symplectic structure on $\bC^2 \times (\bC^\times)^2$ to $X$ because the constructions in \cite{AAK12} are more adapted to this blowup symplectic structure.}

As in the 2-dimensional case, we have a Hamiltonian $S^1$-action
\begin{equation*}
e^{\mathbf{i} \theta} \cdot (x, y, z_1, z_2) = (e^{\mathbf{i}\theta}x, e^{-\mathbf{i}\theta}y, z_1, z_2)
\end{equation*}
and the associated moment map $\mu: X \to \bR_{>0}$.
For each $\lambda \in \bR_{>0}$, Abouzaid-Auroux-Katzarkov \cite{AAK12} showed that there exists a symplectomorphism\footnote{The map $\phi_\epsilon$ is a homoeomorphism which is a diffeomorphism only away from $Z$; see \cite[Lemma 4.1]{AAK12}.}
$$\phi_\lambda: X_{\text{red},\lambda} \to (\bC^\times)^2$$
intertwining between the reduced symplectic form on the reduced space
$$X_{\text{red},\lambda} := \mu^{-1}(\lambda)/S^1$$
and (a constant multiple of) the standard symplectic structure on $(\bC^\times)^2$. Now the map
\begin{equation*}
\rho: X \to B := \bR^2 \times \bR_{>0}, \quad p\in \mu^{-1}(\lambda) \mapsto \left((\text{Log}\circ\phi_\lambda)([p]), \lambda\right),
\end{equation*}
where $\text{Log}: (\bC^\times)^2 \to \bR^2$ is the usual log map,
defines a Lagrangian torus fibration on $X$ \cite[Section 4]{AAK12}.

The discriminant locus of the fibration $\rho$ is given by the {\em amoeba-shaped} subset
$$\Gamma = \mathcal{A} \times \{\epsilon\} = \left(\text{Log}\circ\phi_{\epsilon}\right)(Z) \times \{\epsilon\}.$$
In contrast to the 2-dimensional case, the discriminant locus here is of real codimension one in $B$,
and the walls, i.e. the locus of the Lagrangian torus fibers which bound nontrivial holomorphic disks, is given by
$$H = \mathcal{A} \times \bR,$$
which is an {\em open subset} (though we still call them walls).
The complement of the walls gives a natural chamber structure, which is invariant under vertical translations, in $B$.

\begin{figure}
\begin{tikzpicture}
\draw[dashed] (4,4) -- (14,2.5);
\draw[dashed] (4,4) -- (10,6);
\draw[thick] (4.1,4) -- (3.9,4) node[left] {$\epsilon$};
\draw[->] (3,0.15) -- (14,-1.5);
\draw[->] (3.1,-0.3) -- (10,2);
\draw[->] (4,-1) -- (4,6.5);
\draw [dashed, ultra thick, cyan] (7,4) to [out=7, in=170] (14,3.8);
\draw [dashed, ultra thick, cyan] (6.85,4.15) to [out=8, in=210] (12.2,5.5);
\draw [dashed, ultra thick, cyan] (12.35,5.48) to [out=211, in=168] (14.2,4);
\path [fill=magenta] (6.85,4.15) to [out=8, in=210] (12.2,5.5) -- (12.35,5.48) to [out=211, in=168] (14.2,4) -- (14, 3.8) to [out=170, in=7] (7,4) -- (6.85,4.15);
\draw[dashed, thick, cyan] (7,-0.5) -- (7,4.5);
\draw[dashed, thick, cyan] (14,4.3) -- (14,-0.7);
\draw[dashed, thick, cyan] (6.85,4.65) -- (6.85,-0.35);
\draw[dashed, thick, cyan] (12.2,6) -- (12.2,1);
\draw[dashed, thick, cyan] (12.35,5.98) -- (12.35,0.98);
\draw[dashed, thick, cyan] (14.2,4.5) -- (14.2,-0.5);
\draw [dashed, ultra thick, cyan] (7,0) to [out=7, in=170] (14,-0.2);
\draw [dashed, ultra thick, cyan] (6.85,0.15) to [out=8, in=210] (12.2,1.5);
\draw [dashed, ultra thick, cyan] (12.35,1.48) to [out=211, in=168] (14.2,0);
\path [fill=cyan] (6.85,0.15) to [out=8, in=210] (12.2,1.5) -- (12.35,1.48) to [out=211, in=168] (14.2,0) -- (14, -0.2) to [out=170, in=7] (7,0) -- (6.85,0.15);
\end{tikzpicture}
\caption{The base of SYZ fibration for $X$}\label{fig:3d_SYZbase}
\end{figure}
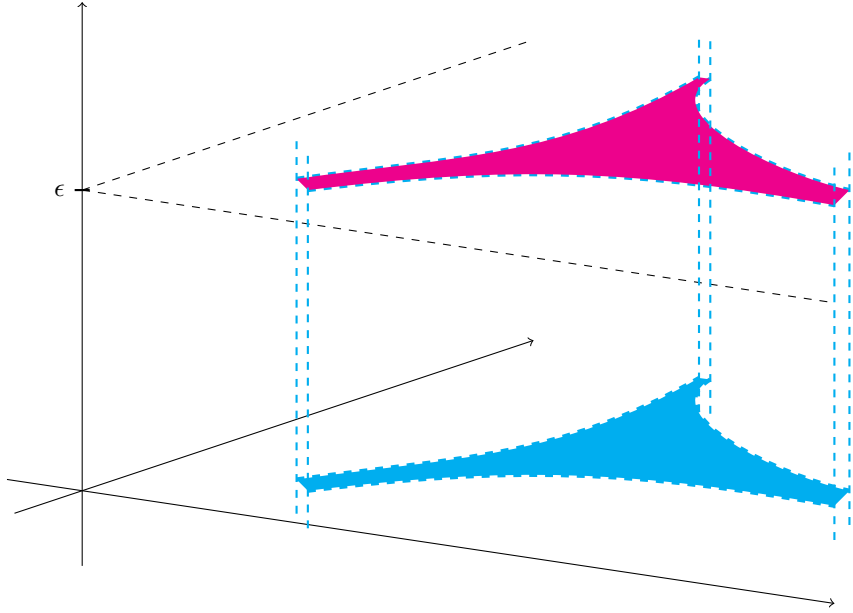

\begin{figure}
\begin{tikzpicture}
\draw [dashed, ultra thick, cyan] (0,4.1) to [out=0, in=225] (5.9,6.1);
\draw [dashed, ultra thick, cyan] (4.1,0) to [out=90, in=225] (6.1,5.9);
\draw [dashed, ultra thick, cyan] (0,3.9) to [out=0, in=90] (3.9,0);
\path [fill=magenta] (0,3.9) to [out=0, in=90] (3.9,0) -- (4.1,0) to [out=90, in=225] (6.1,5.9) -- (5.9,6.1 ) to [out=225, in=0] (0,4.1) -- (0,3.9);
\draw [dashed, very thick] (0,4) -- (4,4);
\draw [dashed, very thick] (4,0) -- (4,4);
\draw [dashed, very thick] (4,4) -- (6,6);
\draw (2,2) node[below] {$U_\alpha$};
\draw (3,5) node[above] {$U_\gamma$};
\draw (5.5,2.8) node[above] {$U_\beta$};
\draw [dashed, very thick] (7,3.5) -- (9,3.5);
\draw [dashed, very thick] (9,1.5) -- (9,3.5);
\draw [dashed, very thick] (9,3.5) -- (10.5,5);
\draw [->] (8.9,3.7) -- (9.4,4.2) node[pos=1, above] {$u_{\gamma,0}$};
\draw [->] (8.9,3.7) -- (8.1,3.7) node[pos=1, above] {$u_{\gamma,1}$};
\draw [->] (8.8,3.3) -- (8,3.3) node[pos=1, below] {$u_{\alpha,1}$};
\draw [->] (8.8,3.3) -- (8.8,2.5) node[pos=1, left] {$u_{\alpha,2}$};
\draw [->] (9.2,3.4) -- (9.7,3.9) node[pos=0.9, right] {$u_{\beta,0}$};
\draw [->] (9.2,3.4) -- (9.2,2.6) node[pos=1, right] {$u_{\beta,2}$};
\end{tikzpicture}
\caption{(Left) Chamber structure in the base of SYZ fibration for $X$.
(Right) Local affine coordinates on $B$.}\label{fig:3d_complement_wall}
\end{figure}
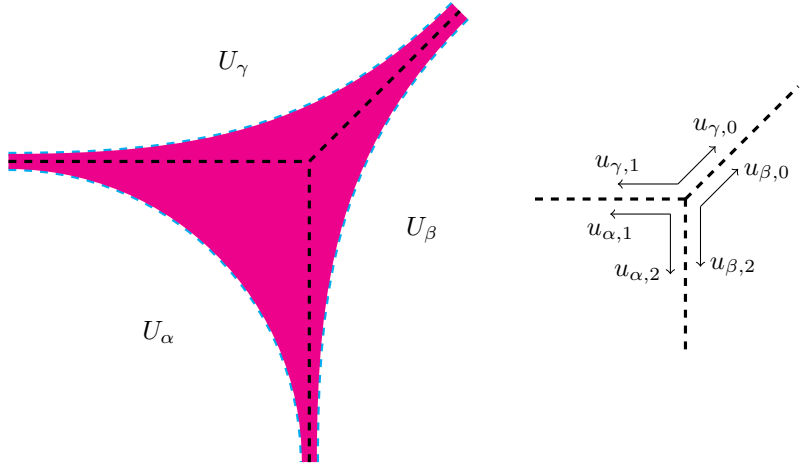

For illustration, let us look at the simplest possible example:
\begin{equation*}
X = \left\{ (x, y, z_1, z_2) \in \bC^2 \times (\bC^\times)^2 \mid xy = 1 + z_1 + z_2 \right\},
\end{equation*}
for which the base of the SYZ fibration is shown in Figure \ref{fig:3d_SYZbase}.
There are 3 chambers in the complement of the walls in $B$ which we denote by
$$U_\alpha, U_\beta, U_\gamma$$
as shown on the left of Figure \ref{fig:3d_complement_wall}.
The semi-flat SYZ mirror $\check{X}_0$ is given by
$$\left( TU_\alpha/\Lambda \right) \cup \left( TU_\beta/\Lambda \right) \cup \left( TU_\gamma/\Lambda \right)$$
using the na\"ive gluing. Note that each chart
$U_{\sharp}^\bC := TU_{\sharp}/\Lambda$ ($\sharp \in \{\alpha, \beta, \gamma\}$) is a copy of the algebraic torus $(\bC^\times)^3$.
Let $u_{\sharp,j}$ ($\sharp \in \{\alpha, \beta, \gamma\}$, $j = 1, 2$) be exponentiation of complexification of the affine coordinates as shown in Figure \ref{fig:3d_complement_wall} (right), and $w$ be exponentiation of complexification of the vertical affine coordinate; again $w$ is a global coordinate since it comes from the moment map of the global Hamiltonian $S^1$-action on $X$.

By analyzing counting of Maslov index two disks (in the blowup $\text{Bl}_{0 \times Z}\left(\bC \times (\bC^\times)^2\right)$, not $X$!), one can see that the corrected gluing should be given by
\begin{equation}\label{eqn:corrected_gluing_2}
\begin{split}
u_{\alpha, 1} & = (1 + w)^{\beta_1 - \alpha_1} u_{\beta, 1}, \\
u_{\alpha, 2} & = (1 + w)^{\beta_2 - \alpha_2} u_{\beta, 2}
\end{split}
\end{equation}
over $U_\alpha^\bC \cup U_\beta^\bC$ and similarly for other pairs of open charts. In the general case, this gives exactly the gluing of charts in the toric Calabi-Yau 3-fold $Y_\Sigma$ whose fan consists of cones over a regular subdivision of the Newton polytope of $f$. Indeed, the above coordinates are related to the coordinates $(t_1, t_2, t_3)$ of the open dense orbit $(\bC^\times)^3 \subset Y_\Sigma$ by
\begin{align*}
u_{\alpha, 1} & = t_1 t_3^{-\alpha_1}, \\
u_{\alpha, 2} & = t_2 t_3^{-\alpha_2}, \\
w & = t_3 - 1,
\end{align*}
and $H$ is the hypersurface defined by $t_3 = 1$.
For the above example, the SYZ mirror precisely is given by
$$\check{X} = \bC^3 \setminus \{w_1w_2w_3 = 1\},$$
where $w_1, w_2, w_3$ are standard coordinates on $(\bC)^3$.

\begin{figure}
\begin{tikzpicture}
\draw [dashed, ultra thick, cyan] (3,4) to [out=0, in=225] (4.5,4.5);
\draw [dashed, ultra thick, cyan] (4,3) to [out=90, in=225] (4.5,4.5);
\draw [dashed, ultra thick, cyan] (3,4) to [out=0, in=90] (4,3);
\path [fill=magenta] (3,4) to [out=0, in=90] (4,3) -- (4,3) to [out=90, in=225] (4.5,4.5) -- (4.5,4.5) to [out=225, in=0] (3,4) -- (3,4);
\draw [dashed, very thick, black] (1,4) -- (4,4);
\draw [dashed, very thick, black] (4,1) -- (4,4);
\draw [dashed, very thick, black] (4,4) -- (6,6);
\draw (3,3) node[below] {$U_\alpha$};
\draw (3,5) node[above] {$U_\gamma$};
\draw (5.5,2.8) node[above] {$U_\beta$};
\draw [ultra thick] (9,3) -- (11,3);
\draw [ultra thick] (9,3) -- (9,5);
\draw [ultra thick] (9,5) -- (11,3);
\draw [dashed, very thick, black] (9.7,3.7) -- (7.7,3.7);
\draw [dashed, very thick, black] (9.7,3.7) -- (9.7,1.7);
\draw [dashed, very thick, black] (9.7,3.7) -- (11,5);
\draw [fill] (9,3) circle [radius=0.1] node[below left] {$\alpha$};
\draw [fill] (11,3) circle [radius=0.1] node[below right] {$\beta$};
\draw [fill] (9,5) circle [radius=0.1] node[above left] {$\gamma$};
\end{tikzpicture}
\caption{(Left) The discriminant locus after tropical localization. (Right) A simplex in the regular subdivision of the Newton polytope of $f$.}\label{fig:3d_trop_localization}
\end{figure}
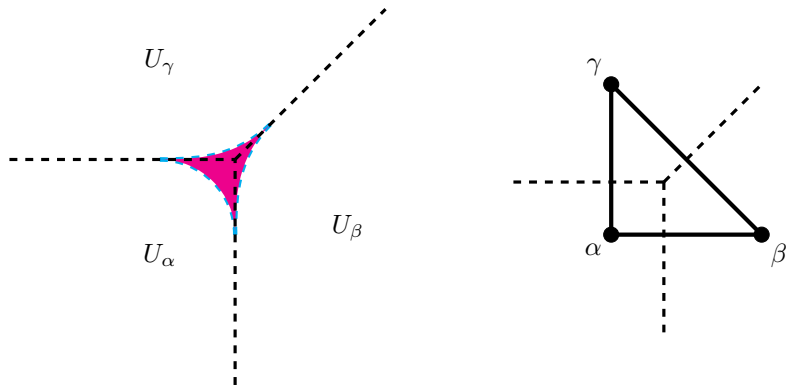

Before describing the SYZ transform in this case, let us briefly explain the class of Lagrangian sections that we are going to look at. We fix a simplex in the regular subdivision of the Newton polytope of $f$ and label its 3 vertices as $\alpha$, $\beta$ and $\gamma$ as shown in Figure \ref{fig:3d_trop_localization} (right).

In \cite{CPU16a}, we employ the technique of so-called {\em tropical localization} devised by Abouzaid in his work \cite{Abouzaid06, Abouzaid09} on HMS for toric varieties. Namely, we consider a family $f_t(z_1, z_2)$ of functions deforming $f(z_1, z_2)$ so that in the limit, where $Z$ becomes $Z^{\text{trop}}$, the ``legs'' of the discriminant locus $\Gamma$ get slimmed down, giving rise to the picture shown in Figure \ref{fig:3d_trop_localization} (left).

We require that our Lagrangian section $L$ is fibered over a Lagrangian submanifold $\underline{L}$ in $(\bC^\times)^2 \setminus U_{Z^{\text{trop}}}$ (where $U_{Z^{\text{trop}}}$ is a tubular neighborhood of $Z^{\text{trop}}$ in $(\bC^\times)^2$), i.e.
\begin{equation*}
L = \bR_{>0} \times \underline{L} \subset \bC \times \left( (\bC^\times)^2 \setminus U_{Z^{\text{trop}}} \right).
\end{equation*}
Note that $p$ is an isomorphism outside of a neighborhood of $0 \times Z^{\text{trop}}$. In fact $\underline{L}$ is a section of the log map $\text{Log}: (\bC^\times)^2 \to \bR^2$ which is disjoint from $U_{Z^{\text{trop}}}$.

We further require that over each chamber $U_\sharp$, $\underline{L}$ is a {\em tropical Lagrangian section}, meaning that its boundary lies in $Z^{\text{trop}}$. In particular, this implies that over a leg $l$ of $\Gamma$ which is dual to the edge connecting two vertices, say $\alpha$ and $\beta$, we have
\begin{equation}\label{eqn:condition_tropical_Lag}
\langle dg, \beta - \alpha \rangle \in \bZ,
\end{equation}
where $g$ denotes the primitive function of $\underline{L}$. This condition is crucial to the construction of the SYZ transform of $L$.

A prototypical example of such a Lagrangian section is given by the positive real locus $L_0 = (\bR_{>0})^3$, which we regard as the zero section of the SYZ fibration $\rho:X \to B$.

As in the 2-dimensional case, over each chamber $U_\sharp$, the semi-flat SYZ transform $\mathcal{F}^\text{semi-flat}$ produces a $U(1)$-connection $\check{\nabla}_\sharp$ over the open chart $U_{\sharp}^\bC = TU_\sharp/\Lambda$ in the mirror $\check{X}$. Due to the modified gluing \eqref{eqn:corrected_gluing_2}, the connections $\check{\nabla}_\sharp$ do not agree over the walls: For example, we have
\begin{equation}\label{eqn:corrected_gauge}
\check{\nabla}_\beta - \check{\nabla}_\alpha = - \mathbf{i} \langle dg, \beta - \alpha \rangle d\text{arg}(1 + w)
\end{equation}
over $l \times \bR$ where $l$ is the leg of $\Gamma$ dual to the edge connecting $\alpha$ and $\beta$.

Remember we have the condition \eqref{eqn:condition_tropical_Lag}: $\langle df, \beta - \alpha \rangle \in \bZ$ over the leg $l$. So, once again, the connections $\check{\nabla}_\sharp$ only differ by a gauge over the intersection of open charts in $\check{X}$. Hence, after modifying by a gauge, they still patch together to give us $U(1)$-connections $\tilde{\nabla}_{\alpha\beta}$, $\tilde{\nabla}_{\beta\gamma}$, $\tilde{\nabla}_{\gamma\alpha}$ over the open charts
\begin{equation*}
U_{\alpha\beta}^\bC = U_\alpha^\bC \cup U_\beta^\bC,\quad U_{\beta\gamma}^\bC = U_\beta^\bC \cup U_\gamma^\bC,\quad U_{\gamma\alpha}^\bC = U_\gamma^\bC \cup U_\alpha^\bC
\end{equation*}
respectively.
Curvatures of these connections have vanishing $(0,2)$-parts as $L$ is Lagrangian, and they satisfy the {\em cocycle condition} since the gauge is of the form specified as in \eqref{eqn:corrected_gauge}. Thus they determine a holomorphic line bundle $\mathcal{L}$ over $\check{X}$, which we define to be the SYZ transform of $L$.

As an example, the SYZ transform of the zero section $L_0 = (\bR_{>0})^3$ is given by the structure sheaf $\mathcal{O}_{\check{X}}$ over $\check{X}$. In \cite{CPU16a}, we computed the wrapped Floer cohomology of $L_0$ and showed that it is isomorphic to the group $H^0(\check{X}, \mathcal{O}_{\check{X}})$ of holomorphic functions over $\check{X}$ (see \cite[Theorem 1.2]{CPU16a}). This result has several interesting applications. We applied it to give a computation of the $0$-th symplectic cohomology of $X$ and a proof of HMS for toric Calabi-Yau orbifold quotients of the form $\bC^3/G$ in \cite[Section 8]{CPU16a}. Recently, it was also used by Ganatra-Pomerleano \cite{Ganatra-Pomerleano16} in giving a complete classification of diffeomorphism types of exact Lagrangian submanifolds in 3-dimensional conic bundles over $(\bC^\times)^2$.

Finally, let us give a couple of remarks regarding the constructions in this article. First of all, we do not expect that the SYZ transform of a general Lagrangian submanifold will be defined in the same way as described here. Indeed one should also take into account the instanton corrections coming from nontrivial holomorphic disks bounded by the Lagrangian itself. In all the examples we consider here (and in the papers \cite{Chan13, Chan-Ueda13, CPU16, CPU16a}), the mirrors of the Lagrangian submanifolds only receive corrections coming from nontrivial holomorphic disks bounded by the SYZ fibers just because the Lagrangians themselves do not bound nontrivial disks. However, this is certainly not the case in general. For instance, in a joint work \cite{Chan-Chung16} with S.-W. Chung, we see that the SYZ transforms of Aganagic-Vafa A-branes \cite{Aganagic-Vafa00} must take into account those disks bounded by the Lagrangian itself.

Related to this, we should point out that a more natural approach is to make direct use of the Floer theory of the Lagrangian submanifolds concerned. More precisely, one may attempt to construct the mirror B-branes (as coherent sheaves over subvarieties) by applying {\em family Floer theory}, as done (in the case without quantum corrections) earlier by Fukaya \cite{Fukaya02b, Fukaya02c} and more recently by Tu \cite{Tu15} and Abouzaid \cite{Abouzaid14a, Abouzaid14b}. Unfortunately, as it is still not known how to deal with family Floer theory for SYZ fibrations containing singular fibers, a general definition of SYZ transforms still eludes us.\footnote{See, however, the recent works of Cho-Lau-Hong \cite{CHL13, CHKL14, CHL14, CHL15} for a generalization of the SYZ approach which hopefully works in more general cases.}

\bibliographystyle{amsplain}
\bibliography{geometry}

\end{document}